\def \CC {\mathbb C}

\def \ZZ {\mathbb Z}

\def \epsilon{\varepsilon}

\def \D  {{\mathcal D}}

\def \S  {{\mathcal S}}

\def \ga {\gamma}

\def \si {\sigma}
\def \la {\lambda}

\newcommand{\res}{{\rm res}}
\renewcommand{\S}{{\mathcal S}}

\documentclass[12pt,reqno]{amsart}

\usepackage{amsmath,amssymb}

\newtheorem{lem}{Lemma}
\usepackage{url} 

\numberwithin{equation}{section}

\begin{document}
\title[]{On the invariants of $L$-functions of degree $2$, I: \\ twisted degree and internal shift}
\author[]{J.Kaczorowski \lowercase{and} A.Perelli}
\maketitle

\medskip
{\bf Abstract.} This is the first part of a series of papers where the behaviour of the invariants under twist by Dirichlet characters  is studied for $L$-functions of degree 2. Here we show, under suitable conditions, that degree and internal shift remain unchanged under twist. The ultimate goal of the series is to prove a general version of Weil converse theorem with minimal assumptions on the shape of the functional equation of the twists.



%-1-%%%%%%%%%%%%%%%%%%%%%%%%%%%%%%%%%%%%%%%%%%%%%%%%%%%%%
%%%%%%%%%%%%%%%%%%%%%%%%%%%%%%%%%%%%%%%%%%%%%%%%%%%%%%%
\section{Introduction}

\smallskip
%%%%%
Assume that 
\begin{equation}
\label{intro-1}
F(s)=\sum_{n=1}^{\infty}a(n)n^{-s}
\end{equation}
belongs to the Selberg class $\S$ and let $\chi(\bmod \, q)$ be a Dirichlet character. In many particular cases, the twist
\[
F^{\chi}(s)=\sum_{n=1}^{\infty}\frac{a(n)\chi(n)}{n^s}
\]
belongs to $\S$ as well. This is true for automorphic $L$-functions modulo some natural restrictions. However, there are no general theorems on $F^\chi$ related to classes of $L$-functions like $\S$ or the extended Selberg class $\S^\sharp$. Several natural questions can be asked on this topic. A general one is the following: Given two $L$-functions $F$ and $F^{\chi}$ belonging ${\mathcal S}^{\sharp}$, what can be said about their invariants like degree, conductor, internal shift, or root number? Can the invariants of the twisted $L$-function be computed in terms of the invariants of the original $L$-function and the data of the character $\chi$?  Based on the known examples, one can risk some conjectures. For instance, one can predict that twisting by a Dirichlet character cannot change the degree. This paper confirms this expectation for a large class of Dirichlet characters and degree two $L$-functions.  

\smallskip
%%%%%
Questions like these are attractive, and are essential for proving converse theorems. For instance, the well-known converse theorem of Weil requires that the twists satisfy very precisely described functional equations, closely related to the functional equation of the untwisted function; cf. \cite[Theorem 7.8]{I-TCAF}. In particular, the invariants of $F^\chi$ can attain determined values only. The present paper is intended as the first part of a series with the aim of proving a general version of the Weil converse theorem with possibly minimal assumptions on the shape of the functional equation of twists. This was the main reason for undertaking the present research.

\smallskip
%%%%%%
Before formulating our main result we introduce the following notation. For a Dirichlet character $\chi$, we denote by $\chi^*$ (mod $f_\chi$) the primitive character inducing $\chi$.  Let $D$ be a fixed positive integer. We say that $F$ is {\it saturated $(\bmod \, D)$} if for every integer $M$ and every reduced residue class $a$ (mod $D$) there exists an integer $\nu$ such that $\nu\equiv a$ (mod $D$), $(\nu,M)=1$ and $a(\nu)\neq 0$. Moreover, we say that $F$ {\it splits polynomially} at a prime $p$ if 
\[ 
F(s)= \prod_{j=1}^{\partial_p}\left(1 - \frac{\alpha_j(p)}{p^s}\right)^{-1} \sum_{p\nmid n}\frac{a(n)}{n^s}.
\]
for some $\partial_p\geq1$ and $\alpha_j(p)\in\CC$. We further recall that the extended Selberg class $\S^\sharp$ consists of the non identically vanishing Dirichlet series $F(s)$ as in \eqref{intro-1}, absolutely convergent for $\si>1$ and satisfying a functional equation of type
\[
\gamma(s) F(s) = \omega \overline{\gamma}(1-s)\overline{F}(1-s),
\]
where $|\omega|=1$ and the $\ga$-factor
\[
 \gamma_F(s)=Q^s\prod_{j=1}^r\Gamma(\lambda_js+\mu_j)
 \]
has $Q>0$, $\lambda_j>0$ and $\Re(\mu_j)\geq 0$ for $j=1,\ldots,r$. The Selberg class $\S$ is the subclass of $\S^\sharp$ of the functions with a general Euler product and satisfying the Ramanujan conjecture. Degree, conductor and internal shift of $F\in\S^\sharp$ are defined respectively as
\[
d_F:=2\sum_{j=1}^r\lambda_j, \qquad q_F:= (2\pi)^{d_F}Q^2\prod_{j=1}^r\lambda_j^{2\lambda_j}, \qquad
\theta_F:= \frac{2}{d_F}\Im\sum_{j=1}^r \mu_j
\]
and are invariants. We refer to our survey \cite{Ka-Pe/1999b} and to the papers by Selberg \cite{Sel/1989}, Conrey-Ghosh \cite{Co-Gh/1993} and Murty \cite{Mur/1994} for further information on the Selberg classes $\S$ and $\S^\sharp$.

\medskip
%%%%%
{\bf Theorem.} {\sl Let $F\in {\mathcal S}^{\sharp}$ have degree $2$, integer conductor $q_F$ and internal shift $\theta_F$. Let $D$ be a positive square-free integer with $D\equiv 1$ {\rm (mod $q_F$)} such that $F$ splits polynomially at every prime $p|D$. Moreover, suppose that $F$ is saturated $(\bmod \, D)$ and for every Dirichlet character $\chi$ {\rm (mod $D$)}, the twist $F^{\chi^*}$ belongs to ${\mathcal S}^{\sharp}$ and has degree $d_{\chi^*}$ and internal shift $\theta_{\chi^*}$. Then $d_{\chi^*}=2$ and $\theta_{\chi^*}=\theta_F$.}  

\medskip
%%%%
We remark that the assumption that $F$ is saturated (mod $D$) can most probably be relaxed, at least in the case of $L$-functions in the Selberg class $\S$. Using the linear independence of Dirichlet characters, one quickly shows that $F$ is saturated (mod $D$) if and only if  the twists $F_{|M}^{\chi}$, where 
\[
F_{|M}(s)= \sum_{(n,M)=1}\frac{a(n)}{n^s}
\]
and $\chi$ (mod $D$), are linearly independent for every positive integer $M$.  Thus, if $F$ and all the twists $F^{\chi^*}$, $\chi$ (mod $D$), belong to $\S$,  then the following statements are equivalent:

(i) $F$ is saturated $(\bmod \, D)$;

(ii) the map $\chi$ (mod $D$) $\mapsto F^{\chi^*}$ is injective, i.e., for  $\chi$ (mod $D$)  all $F^{\chi^*}$ are different.

\noindent
This follows from a general fact, namely that distinct functions from the Selberg class $\S$ are linearly independent over the ring of $p$-finite Dirichlet series; cf. \cite[Theorem 1]{Ka-Mo-Pe/1999},. Actually, we expect all functions in $\mathcal S$ to be saturated $(\bmod \, D)$ for every positive integer $D$.

\medskip
%%%%%%
{\bf Acknowledgements}.  We wish to thank the referee for the suggestions leading to a better presentation of the paper. This research was partially supported by the Istituto Nazionale di Alta Matematica, by the MIUR grant PRIN-2017 {\sl ``Geometric, algebraic and analytic methods in arithmetic''} and by grant 2021/41/BST1/00241 {\sl ``Analytic methods in number theory''}  from the National Science Centre, Poland.

\section{Notation.} 
%-2-%%%%%%%%%%%%%%%%%%%%%%%%%%%%%%%%%%%%%%%%%%%%%%%%%%%%%%%

For real numbers $\alpha$, $\beta$, $\lambda$ and $F\in {\mathcal S}^{\sharp}$ we consider the nonlinear twist
\begin{equation}
\label{nonlindef}
F(s,\beta,\alpha,\lambda) = \sum_{n=1}^{\infty}\frac{a(n)}{n^s}e(-\beta n -\alpha n^{\lambda}).
\end{equation}
If $d_F>0$ and $\alpha>0$, then $F(s,0,\alpha,1/d_F)$ becomes the standard twist of $F$, is denoted by $F(s,\alpha)$ and its properties are described in Theorems 1 and 2 in \cite{Ka-Pe/2005}. Indeed, $F(s,\alpha)$ admits meromorphic continuation to the whole complex plane, and in the half-plane $\sigma>1/2-1/(2d_F)$, it has at most a simple pole at $s=s_0:=1/2+1/(2d_F)-i\theta_F$ with residue
\begin{equation}\label{eq:residue}
{\res}_{s=s_0} F(s,\alpha) =c(F)\frac{\overline{a(n_{\alpha})}}{n_{\alpha}^{1-s_0}},
\end{equation}
where $c(F)\neq 0$ is a constant depending on $F$ and 
\begin{equation}\label{eq:nalpha}
n_\alpha = q_Fd_F^{-d_F}\alpha^{d_F}.
\end{equation} 
We use the convention that $a(\xi)=0$ if $\xi\not\in {\mathbb N}$. 
 In particular, $F(s,\alpha)$ is holomorphic at $s_0$ whenever $n_{\alpha}$ is not a positive integer. We also write $\overline{F}(s) = \overline{F(\overline{s})}$.

\smallskip
%%%%%
For a prime $p$ we denote by $F_p(s)$ the corresponding local factor of $F$, i.e. $F_p(s)=\sum_{\ell=0}^{\infty} a(p^{\ell})p^{-\ell s}$.
So $F(s)$ splits polynomially at a prime $p$ if $F(s)= F_{|p}(s) F_p(s)$ and there exist an integer $\partial_p\geq 0$ and coefficients 
\begin{equation}\label{eq:Ap}
A_{\ell}(p)=A_{\ell}(p, F), \qquad \ell=0,\ldots,\partial_p \ \text{and} \ A_{\partial_p}(p)\neq 0
\end{equation}
 such that
\[ 
F_p^{-1}(s)=\sum_{\ell=0}^{\partial_p} A_{\ell}(p) p^{-\ell s}.
\]
If $F$ splits polynomially at all prime factors of a positive integer $D$ we write
\[ D^*=D^*(F)=\prod_{p|D}p^{\partial_p}\]
and
\begin{equation}\label{eq:Bm}
 \prod_{p|D}F_p^{-1}(s)=\sum_{m|D^*} B_m(F) m^{-s},
 \end{equation}
where
\begin{equation}\label{eq:Bm'}
B_m(F)=\prod_{p^{\ell}||m} A_{\ell}(p).
\end{equation}
Clearly, $B_m(F)$ is multiplicative in $m$.

\smallskip
%%%%%
Finally, $\tau(\chi^*)$ denotes the Gauss sum of $\chi^*$, and the radical $r(m)$ of an integer $m=p_1^{a_1}\cdots p_k^{a_k}$ is defined as
\[
r(m) = p_1\cdots p_k.
\]

%-3-%%%%%%%%%%%%%%%%%%%%%%%%%%%%%%%%%%%%%%%%%%%%%%%%%%%%%%%
\section{Lemmas}

%%%%%%%%%%%%%%%%%%%%%%%
\begin{lem}\label{lem:chi0e}
Let $D$ be a positive integer, $\chi_0$ be the principal character {\rm (mod $D$)} and let $(a,D)=1$. Then for every positive integer $n$ we have
\begin{equation}\label{eq:chi0e}
\chi_0(n)e(-an/D) =\frac{1}{\varphi(D)}\sum_{\chi(\bmod \, D)} c(\chi, a/D) \chi(n),
\end{equation}
where
\begin{equation}\label{eq:cchi}
 c(\chi,a/D)= \mu(D/f_{\chi}) \chi^*(a) \overline{\chi^*(D/f_{\chi})} \,  \overline{\tau(\chi^*)}.
 \end{equation}
\end{lem}
%%%%%%%%%%%%%%%%%%%%%%%

\smallskip
%%%%
{\it Proof.} By elementary Fourier analysis (\ref{eq:chi0e}) holds with 
\[
c(\chi,a/D) =\sum_{n(\bmod \, D)} \overline{\chi(n)}e(-an/D).
\] 
But with the notation in Chapter 9 of Montgomery-Vaughan \cite{M-V} we have that $c(\chi,a/D)=\overline{c_\chi(a)}$, hence \eqref{eq:cchi} follows at once from Theorem 9.12 in \cite{M-V}. \qed

\smallskip
%%%%%%%%%%%%%%%%%%%%%%%
\begin{lem}\label{lem:chi0} Let $D$ be a positive squarefree integer, suppose that $F\in \mathcal{S}^{\sharp}$ splits polynomially at all primes $p|D$ and let $\chi_0$ be the principal character {\rm (mod $D$)}. Then for every positive integer $n$ we have
\[
a(n) \chi_0(n) = \sum_{m|(n, D^*)} B_m(F) a(n/m).
\]
\end{lem}
%%%%%%%%%%%%%%%%%%%%%%%

\smallskip
%%%%
{\it Proof.} We have 
\[ 
\sum_{n=1}^{\infty}\frac{a(n) \chi_0(n)}{n^s} = \left(\prod_{p|D} F^{-1}_p(s)\right)F(s),
\]
and the lemma follows at once recalling the expression (\ref{eq:Bm}) for the function inside brackets. \qed

\smallskip
%%%%%%%%%%%%%%%%%%%%%%%
\begin{lem}\label{lem:fchi}
Let $D$ be a positive squarefree integer and $(a,D)=1$, and suppose that $F\in \mathcal{S}^{\sharp}$ splits polynomially at all primes $p|D$. Then for $\si>1$ we have
\begin{equation}\label{eq:fchi1}
F(s,a/D,\alpha,\lambda) = \frac{1}{\varphi(D)} \sum_{\chi(\bmod \, D)} \sum_{m|\left(\frac{D}{f_{\chi}}\right)^*}
\frac{B_m(F)f(\chi,m,a)}{m^s} F^{\chi^*}(s,0,m^{\lambda}\alpha,\lambda),
\end{equation}
where
\begin{equation}\label{eq:fchi2}
f(\chi,m,a)= \mu\left(\frac{D}{f_{\chi}}\right) \overline{\tau(\chi^*)}\overline{\chi^*\left(\frac{D}{f_{\chi}}\right)}
\chi^*(am) r(m) \neq 0
\end{equation}
and $\chi^*\,  (${\rm mod} $f_\chi)$ is the primitive Dirichlet character inducing the character $\chi$.
\end{lem}
%%%%%%%%%%%%%%%%%%%%%%%

\smallskip
%%%%%
{\it Proof.}  We proceed by induction over $D$. If $D=1$ the assertion is obvious since both sides of (\ref{eq:fchi1}) equal $F(s,0,\alpha,\lambda)$. Indeed, in this case we have $B_1(F)=1$, $f(\chi,1,1)=1$ and $\chi^*(n)=1$ for every $n$.

\smallskip
%%%%%
Let now $D>1$ and assume the assertion for all proper divisors of $D$. By Lemma \ref{lem:chi0} we have
\begin{eqnarray*}
F_{|D}(s, a/D,\alpha,\lambda) &=&  
\sum_{n=1}^{\infty} \frac{a(n)\chi_0(n)}{n^s} e(-{an}/{D}-\alpha n^{\lambda})\\
&=&\sum_{n=1}^{\infty}n^{-s}\sum_{k|(n,D^*)} B_k(F)a({n}/{k}) e(-{an}/{D}-\alpha n^{\lambda})\\
&=& \sum_{k|D^*} \frac{B_k(F)}{k^s} F(s, {ak}/{D}, k^{\lambda}\alpha, \lambda),
\end{eqnarray*}
Hence
\begin{equation}\label{eq:A-B}
\begin{split}
F(s,a/D,\alpha,\lambda) &= F_{|D}(s, a/D,\alpha,\lambda)
-\sum_{\substack{k|D^*\\ k>1}}\frac{B_k(F)}{k^s} F(s, {ak}/{D}, k^{\lambda}\alpha, \lambda) = A - B.
\end{split}
\end{equation}

\smallskip
%%%%%
We consider $A$ first. By Lemma \ref{lem:chi0e} we have
\[
A= \sum_{n=1}^\infty \frac{a(n)\chi_0(n)}{n^s}e(-an/D)e(-\alpha n^\lambda) = \frac{1}{\varphi(D)}\sum_{\chi(\bmod D)}c(\chi, a/D) F^{\chi}(s, 0,\alpha,\lambda).
\]
Observing that $F^{\chi}(s, 0,\alpha,\lambda)=F_{|D}^{\chi^*}(s, 0,\alpha,\lambda)$, applying  Lemma \ref{lem:chi0} and (\ref{eq:cchi}), we obtain 
\begin{equation}\label{eq:A}
A= \frac{1}{\varphi(D)}\sum_{\chi(\bmod D)} \sum_{m|\left(\frac{D}{f_{\chi}}\right)^*}
\frac{B_m(F) \mu(D/f_{\chi}) \overline{\chi^*(D/f_{\chi})} \, \overline{\tau(\chi^*)}\chi^*(am)}{m^s} F^{\chi^*}
(s,0,m^{\lambda}\alpha,\lambda),
\end{equation}
which is of the form \eqref{eq:fchi1} and \eqref{eq:fchi2} with 1 in place of $r(m)$. Next we consider $B$. Since 
\begin{equation}
\label{inductive}
F(s, ak/D, k^{\lambda}\alpha, \lambda) =  F\left(s, \frac{ak/(k,D)}{D/(k,D)}, k^{\lambda}\alpha, \lambda\right)
\end{equation}
and $k>1$, $k|D^*$, we have that $D/(k,D)$ is a proper divisor of $D$ and $(\frac{ak}{(k,D)}, \frac{D}{(k,D)})=1$, so we may apply the inductive assumption to the right hand side of \eqref{inductive}. We obtain
\begin{equation}\label{eq:B}
\begin{split}
B= \frac{1}{\varphi(D)}&\sum_{\substack{k|D^* \\ k>1}} \varphi((k,D))
\sum_{\chi(\bmod \frac{D}{(k,D)})} \mu(D/f_{\chi})\overline{\chi^*(D/f_{\chi})} \, \overline{\tau(\chi^*)}\chi^*(a)
\mu((k,D))\\
&\times \sum_{\nu|\left(\frac{D}{f_{\chi}(k,D)}\right)^*}\frac{B_{k\nu}(F)\chi^*(k\nu)r(\nu)}{(k\nu)^s}
F^{\chi^*}(s,0, (k\nu )^{\lambda}\alpha,\lambda),
\end{split}
\end{equation}
since $B_m(F)$ is multiplicative and $(k,\nu)=1$. But $D$ is squarefree, so if $k|D^*$, $\nu |(D/(f_{\chi}(k,D))^*$ and $\chi^*(k)\neq 0$ then
\[ 
k\nu \ \text{divides} \  \left(\frac{D}{f_{\chi}}\right)^*.
\]
We also observe that the sets of primitive characters $\chi^*$ obtained from the characters $\chi$ (mod $D$), and from the characters $\chi$ (mod $D/d$), as $d$ ranges over the divisors of $D$, are the same. Thus, writing $m=k\nu$, (\ref{eq:A-B}),\eqref{eq:A} and (\ref{eq:B}) imply (\ref{eq:fchi1}) with certain coefficients $f(\chi,m, a)$. It remains to show that they satisfy (\ref{eq:fchi2}). 

\smallskip
%%%%%
If $f_{\chi}=D$ then $B$ does not contribute any $\chi^*$ to \eqref{eq:fchi1} since $D/(k,D)<D$ in \eqref{eq:B}. Moreover, the inner sum in $A$ has only single term corresponding to  $m=1$. So, in this case 
\[
f(\chi,1,a)=\mu\left(\frac{D}{f_{\chi}}\right)  \overline{\chi^*\left(\frac{D}{f_{\chi}}\right)}\overline{\tau(\chi^*)}
\chi^*(a),
\]
as required. Assume now $f_{\chi}<D$. In this case both $A$ and $B$ contribute to $f(\chi,m, a)$.   From (\ref{eq:A}) we see that $A$ contributes
\begin{equation}\label{eq:fA}
\mu\left(\frac{D}{f_{\chi}}\right) \overline{\chi^*\left(\frac{D}{f_{\chi}}\right)}\overline{\tau(\chi^*)}\chi^*(am),
\end{equation}
while the contribution of $B$ is
\begin{equation}\label{eq:fB}
\mu\left(\frac{D}{f_{\chi}}\right) \overline{\chi^*\left(\frac{D}{f_{\chi}}\right)}\overline{\tau(\chi^*)}\chi^*(am) \sum_{\substack{k|m\\k>1\\ (k,m/k)=1}}\varphi((k,D))\mu((k,D)) r\left(\frac{m}{k}\right).
\end{equation}
But
\[
\begin{split}
r(m) \sum_{\substack{k|m\\ (k,m/k)=1}}\varphi((k,D))\mu((k,D)) r\left(\frac{m}{k}\right)/r(m) &= r(m) \prod_{p^a\|m}\left(1+ \frac{\varphi((p^a,D))\mu((p^a,D))}{p}\right) \\
&= r(m) \prod_{p^a\|m}\left(1-\frac{p-1}{p}\right) = 1,
\end{split}
\]
hence \eqref{eq:fB} becomes
\[
\mu\left(\frac{D}{f_{\chi}}\right) \overline{\chi^*\left(\frac{D}{f_{\chi}}\right)}\overline{\tau(\chi^*)}\chi^*(am)(1-r(m)),
\]
and the lemma follows from \eqref{eq:A-B},\eqref{eq:fA} and \eqref{eq:fB}. \qed

\smallskip
%%%%%%%%%%%%%%%%%%%%%%%
\begin{lem}\label{lem:F-Prop}
Let $F\in {\mathcal S}^{\sharp}$ have degree $2$ and integer conductor $q_F$. Moreover, let $D\equiv 1(\bmod \, q_F)$ be a positive integer. Then the nonlinear twist $F(s, 1/D, \alpha,\lambda)$ is entire if $0<\lambda<1/2$, and is holomorphic for $s\in {\mathbb C}\backslash(-\infty,3/4-i\theta_F]$ if $\lambda =1/2$.
\end{lem}
%%%%%%%%%%%%%%%%%%%%%%%

\smallskip
%%%%%
{\it Proof.} We use the notation in our paper \cite{Ka-Pe/2016b}. Writing 
\begin{equation}
\label{L4}
f(n)=n/D+\alpha n^\lambda, 
\end{equation}
by \cite[Theorem 1]{Ka-Pe/2016b}, we know that the term $F(s;f) = F(s,1/D,\alpha,\lambda)$ can be expressed in terms of the dual twist 
\[
\overline{F}(s;f^*) = \sum_{n=1}^\infty \frac{\overline{a(n)}}{n^s} e(f^*(n))
\]
in the following way. For every $K>0$ there exists an integer $J=J(K)>0$, functions $W_0(s),\dots,W_J(s)$ and $G_J(s;f)$ holomorphic for $\si>-K$, and real numbers $0=\eta_0<\eta_1<\dots<\eta_J$ such that
\begin{equation}
\label{L4-0}
F(s,1/D,\alpha,\lambda)= \sum_{j=0}^J W_j(s) \overline{F}(s+\eta_j+2i\theta_F;f^*) +G_J(s;f).
\end{equation}
Here $\theta_F$ is the internal shift defined in Section 2 and $f^*$ is the dual of $f$. The precise shape of $f^*$ is important in this lemma, but only a general description of $f^*$ is given in \cite{Ka-Pe/2016b}; see Section 1.3 there. So we first proceed with the explicit computation of $f^*$.

\smallskip
%%%%%
We consider, more generally, the function
\begin{equation}
\label{L4-0'}
f(n) = \beta n + \alpha n^\la.
\end{equation}
The first step requires the computation of the $z$-critical point $x_0=x_0(\xi)$ of the function $\Phi(z,\xi)$ defined by equation (1.5) of \cite{Ka-Pe/2016b}, which in our case becomes
\begin{equation}
\label{L4-1}
\Phi(z,\xi) = z^{1/2} - 2\pi \frac{q\beta z}{\xi} -2\pi\alpha \frac{q^\la z^\la}{\xi^\la}, \hskip1cm q=\frac{q_F}{(4\pi)^2}.
\end{equation}
Hence $x_0$ satisfies
\begin{equation}
\label{L4-2}
x_0^{1/2} = 4\pi q\beta \frac{x_0}{\xi} +4\pi\alpha\la q^\la \left(\frac{x_0}{\xi}\right)^\la = 4\pi q\beta \frac{x_0}{\xi} \psi\left(\frac{x_0}{\xi}\right), \hskip1cm \psi(u) = 1+\frac{\alpha\la}{\beta} q^{\la-1} u^{\la-1}.
\end{equation}
Consequently, squaring both sides and solving for $(x_0/\xi)^2$, we get
\begin{equation}
\label{L4-3}
\frac{x_0}{\xi} = (4\pi q\beta)^{-2} \xi \, \psi(x_0/\xi)^{-2}.
\end{equation}

\smallskip
%%%%%
Next we recall from equation (1.8) of \cite{Ka-Pe/2016b} that
\begin{equation}
\label{L4-4}
f^*(\xi) = \frac{1}{2\pi} \Phi^\flat(x_0,\xi),
\end{equation}
where in our case $\Phi^\flat$ means that we drop all terms with $\omega>1$ of the asymptotic expansion
\[
\frac{1}{2\pi}\Phi(x_0,\xi) = \xi\sum_{\omega\in\D_f} A_\omega \xi^{-\omega}, \hskip1cm \D_f = \{\omega=m(1-\la): m\in\ZZ,m\geq0\};
\]
see Section 1.3 of \cite[p.6742-6743]{Ka-Pe/2016b}. Hence, if $0<\la<1/2$ we have
\begin{equation}
\label{L4-4'}
\frac{1}{2\pi} \Phi^\flat(x_0,\xi) = A_0\xi + A_{1-\la} \xi^\la,
\end{equation}
while if $\la=1/2$ there is an additional constant term, which we omit since it is clearly irrelevant in \eqref{L4-0}. Now we compute the coefficients $A_0$ and $A_{1-\la}$.

\smallskip
%%%%%
From \eqref{L4-1},\eqref{L4-2} and \eqref{L4-3} we get
\begin{equation}
\label{L4-5}
\begin{split}
\frac{1}{2\pi}\Phi(x_0,\xi) &= 2q\beta \frac{x_0}{\xi} \psi\left(\frac{x_0}{\xi}\right) -q\beta \frac{x_0}{\xi} \left(1+\frac{\alpha}{\beta} q^{\la-1}\left(\frac{x_0}{\xi}\right)^{\la-1}\right) \\
&= q\beta\frac{x_0}{\xi} \left(1- \frac{\alpha}{\beta}(1-2\la)q^{\la-1}\left(\frac{x_0}{\xi}\right)^{\la-1}\right) \\
&= \frac{1}{q_F\beta} \xi \psi^{-2} \left(\frac{x_0}{\xi}\right) \left(1- \frac{\alpha}{\beta}(1-2\la) q^{\la-1}\left(\frac{x_0}{\xi}\right)^{\la-1}\right).
\end{split}
\end{equation}
Moreover, recalling the definition of $\psi(u)$ in \eqref{L4-2}, as $u\to\infty$ we have
\[
\psi^{-2}(u) = 1- \frac{2\alpha\la}{\beta} q^{\la-1} u^{\la-1} + O(u^{2(\la-1)}),
\]
hence
\begin{equation}
\label{L4-6}
\psi^{-2}(u) \left(1- \frac{\alpha}{\beta}(1-2\la) q^{\la-1} u^{\la-1}\right) = 1- \frac{\alpha}{\beta} q^{\la-1} u^{\la-1} + O(u^{2(\la-1)}).
\end{equation}
But from \cite[Lemma 2.3]{Ka-Pe/2011a} we know that in our case
\begin{equation}
\label{L4-7}
\frac{x_0}{\xi} = \left(\frac{4\pi}{\beta q_F}\right)^2\xi +O(\xi^{1-\eta}), \quad \text{with some $\eta>0$}.
\end{equation}
Thus, recalling the definition of $q$ in \eqref{L4-1}, from \eqref{L4-5},\eqref{L4-6} and \eqref{L4-7} we obtain that
\begin{equation}
\label{L4-8}
\frac{1}{2\pi}\Phi(x_0,\xi) = \frac{1}{q_F\beta} \xi - \frac{\alpha}{\beta^{2\la} q_F^\la} \xi^\la + o(\xi^\la).
\end{equation}
Therefore, from \eqref{L4-4},\eqref{L4-4'} and \eqref{L4-8} we finally get that
\begin{equation}
\label{L4-9}
f^*(n) = \frac{1}{q_F\beta} n - \frac{\alpha}{\beta^{2\la} q_F^\la} n^\la.
\end{equation}

\smallskip
%%%%%
We are now ready to conclude the proof of the lemma. Recalling definition \eqref{nonlindef}, thanks to \eqref{L4},\eqref{L4-0},\eqref{L4-0'} and \eqref{L4-9} we have that
\[ 
\begin{split}
F(s,1/D,\alpha,\lambda) &= \sum_{j=0}^J W_j(s) \overline{F}(s+\eta_j+2i\theta_F, D/q_F, -\alpha D^{2\lambda}q_F^{-\lambda}, \lambda) + G_J(s;f) \\
& = \sum_{j=0}^J W_j(s) \overline{F}(s+\eta_j+2i\theta_F, 1/q_F, -\alpha D^{2\lambda}q_F^{-\lambda}, \lambda) + G_J(s;f)
\end{split}
\]
since $D\equiv 1$ (mod $q_F$).
Then, applying once again transformation formula \eqref{L4-0} to the twists in the last equation, we finally obtain an expression of type
\begin{equation}
\label{L4-10}
F(s,1/D,\alpha,\lambda) = \sum_{j=0}^{\widetilde{J}} \widetilde{W}_j(s) F(s+\widetilde{\eta}_j,0,c, \lambda) + \widetilde{G}_{\widetilde{J}}(s;f)
\end{equation}
with the same properties as \eqref{L4-0}, where $c>0$ is a certain constant. Now, we suppose that $0<\la<1/2$. Then, using the notation of \cite[p.1352-1353]{Ka-Pe/2016a}, thanks to Remark 5 there, we have that $F$ belongs to the class $M(2,\tau)$ for every $\tau>0$. Thus, by \cite[Theorem 5]{Ka-Pe/2016a}, all twists on the right hand side of \eqref{L4-10} are entire. Therefore $F(s,1/D,\alpha,\lambda)$ is entire as well since $K$ is arbitrary, and the first assertion of the lemma follows. The second assertion follows from well known properties of the standard twist, see \cite[Theorems 1 and 2]{Ka-Pe/2016a}. \qed

%-4-%%%%%%%%%%%%%%%%%%%%%%%%%%%%%%%%%%%%%%%%%%%%%%%%%%%%%%%%%%%%%
\section{Proof of the theorem} 

Recall that for $F$ as in the Theorem and a Dirichlet character $\chi(\bmod \, D)$, $d_{\chi^*}$ and $\theta_{\chi^*}$ denote degree and internal shift of the twist $F^{\chi^*}(s)$, respectively.

\smallskip
%%%%%
We first observe that $d_{\chi^*}\geq 2$ for all $\chi(\bmod \, D)$. Indeed, otherwise we would have $d_{\chi^*}=0$ or $d_{\chi^*}=1$, as there are no $L$-functions in ${\mathcal S}^{\sharp}$ of degrees $0<d<1$ and $1<d<2$, see \cite{Ka-Pe/1999a} and \cite{Ka-Pe/2011a}. Both possibilities lead to a contradiction. In the first case, $F^{\chi^*}(s)$ is a Dirichlet polynomial (see \cite[Theorem 1]{Ka-Pe/1999a})
\[ 
F^{\chi^*}(s) = D(s),
\]
and in the second a linear combination of shifted Dirichlet $L$-functions corresponding to the characters (mod $q(\chi^*)$)
\[ 
F^{\chi^*}(s)= \sum_{\psi (\bmod q(\chi^*))} P(s,\psi) L(s+i \theta, \psi^*),
\]
$q(\chi^*)$ being the conductor of $F^{\chi^*}$ (see \cite[Theorem 2]{Ka-Pe/1999a}). Here $P(s,\psi)$ are certain Dirichlet polynomials and $\theta$ is a certain real shift. Thus, twisting by $\overline{\chi^*}$  we have
\begin{equation}\label{eq:d=0}
 F^{\chi_0}(s)= D^{\overline{\chi^*}}(s)
 \end{equation}
and
\begin{equation}\label{eq:d=1}
F^{\chi_0}(s)=\sum_{\psi(\bmod q(\chi^*))} P^{\overline{\chi^*}}(s,\psi) L(s+i \theta, \psi \overline{\chi^*}),
\end{equation}
respectively, where $\chi_0$ denotes the principal character $(\bmod \, f_{\chi})$. Obviously, $f_{\chi}|D$ and hence
\[ 
F^{\chi_0}(s) = F(s) \prod_{p|f_{\chi}} F_p^{-1}(s).
\]
Recalling that the Lindel\"of $\mu$-function of a given function $f(s)$ is
\[
\mu(\si,f) = \inf\{\xi\geq0: f(\si+it)\ll (|t|+1)^\xi\},
\]
we therefore see that the Lindel\"of $\mu$-functions of $F^{\chi_0}(s)$ and $F(s)$ coincide. In particular, for $\sigma<0$ we have
\[ 
\mu(\sigma, F^{\chi_0}) = 2\left(\frac{1}{2}-\sigma\right).
\]
But computing the same function using (\ref{eq:d=0}) and (\ref{eq:d=1}), for $\si<0$ we obtain that
\[
 \mu(\sigma, F^{\chi_0}) \leq \frac{1}{2}-\sigma,
\]
a contradiction. Thus $d_{\chi^*}\geq 2$ for all $\chi(\bmod \, D)$.

\smallskip
%%%%%
The rest of the proof will be conducted by contradiction as well. Suppose that
\[ 
d_0:=\max_{\chi(\bmod D)} d_{\chi^*} >2
\]
and write
\[
A_0=\{ \chi(\bmod D): d_{\chi^*}=d_0\},\quad  \theta_0= \min_{\chi\in A_0} \theta_{\chi^*}, \quad B_0=\{\chi\in A_0: \theta_{\chi^*}=\theta_0\}, \quad \la_0 = 1/d_0<1/2.
\]
By Lemma \ref{lem:F-Prop}, we know that the function $F(s,1/D,\alpha,\lambda_0)$ is entire. On the other hand, according to Lemma \ref{lem:fchi}, it is a combination of twists
 \[ 
 F^{\chi^*}(s,0,m^{\lambda_0}\alpha,\lambda_0)
 \]
with $\chi(\bmod \, D)$ and $m|(D/f_{\chi})^*$. If $\chi\not\in A_0$ then $\lambda_0<1/d_{\chi^*}$, thus the corresponding term is entire thanks to \cite[Theorem 5]{Ka-Pe/2016a}. The other terms are standard twists of $F^{\chi^*}$, so may have a pole at $s_0(\chi^*)=1/2+1/(2d_0) - i\theta_{\chi^*}$, but from (\ref{eq:fchi1}) we have
 \begin{equation}\label{eq:sumres=0}
 \sum_{\chi\in B_0}\sum_{m|\left(\frac{D}{f_{\chi}}\right)^*} \frac{B_m(F)f(\chi,m,1/D)}{m^{s_0}} {\rm Res}_{s=s_0} F^{\chi^*}(s,0,m^{\lambda_0}\alpha,\lambda_0) =0
 \end{equation}
where $s_0=1/2+1/(2d_0) - i\theta_0$. Recalling (\ref{eq:residue}) and (\ref{eq:nalpha}) we have
\begin{equation}\label{eq:residue2}
{\rm Res}_{s=s_0} F^{\chi^*}(s,0,m^{\lambda_0}\alpha,\lambda_0)=
c(F^{\chi^*})\frac{\overline{a(n_{m^{\lambda_0}\alpha})}\,  \overline{\chi^*(n_{m^{\lambda_0}\alpha})}}{n_{m^{\lambda_0}\alpha}^{1-s_0}},
\end{equation}
where $c(F^{\chi^*})\neq 0$ and
\[
n_{m^{\lambda_0}\alpha}=n_{m^{\lambda_0}\alpha}(\chi^*)=q({\chi^*)}d_0^{-d_0}\alpha^{d_0}m.
\]

\smallskip
%%%%%
Let now
\[
q_0= \max_{\chi\in B_0}\,  (D/f_{\chi})^*q({\chi^*}), \ \ \text{and for positive integers $\nu$ let $ \alpha_{\nu}=d_0(\nu/q_0)^{1/d_0}$.}
\]
If $\chi\in B_0$ then 
%%
%\begin{equation}\label{eq:nalpha2}
\[
n_{m^{\lambda_0}\alpha_{\nu}}(\chi^*)= \frac{q(\chi^*)m}{q_0} \nu,
\]
%\end{equation}
%%
hence if $\chi\in B_0$  is such that $q(\chi^*)/q_0\not\in {\mathbb Q}$, then $a(n_{m^{\lambda_0}\alpha_{\nu}})=0$ and the residue in (\ref{eq:residue2}) vanishes. Let further
\[ 
C_0=\{\chi\in B_0: \left(\frac{D}{f_{\chi}}\right)^* q(\chi^*) = q_0\},
\]
$\chi\in B_0\backslash C_0$ and 
\begin{equation}\label{eq:B0-C0}
\frac{q(\chi^*)}{q_0} = \frac{l_{\chi^*}}{m_{\chi^*}}, \qquad l_{\chi^*},m_{\chi^*} \in {\mathbb N}, \, (l_{\chi^*},m_{\chi^*})=1,\, m_{\chi^*}>1.
\end{equation}
Then 
\[
n_{m^{\lambda_0}\alpha_{\nu}}(\chi^*)= \frac{l_{\chi^*}m}{m_{\chi}^*} \nu\not\in {\mathbb N}
\]
for all $(\nu,m_{\chi^*})=1$. Denoting by $M$ the least common multiple of all numbers $m_{\chi^*}$ in (\ref{eq:B0-C0}) we have $a(n_{m^{\lambda_0}\alpha_{\nu}}) =0$ for all $(\nu,M)=1$, and the corresponding residues in (\ref{eq:residue2}) vanish in all these cases. Suppose now that $\chi\in C_0$ and $m$ is a proper divisor of $(D/f_{\chi})^*$. Then for $(\nu, D)=1$ we have
\[
 n_{m^{\lambda_0}\alpha_{\nu}}(\chi^*) = \frac{m\nu}{\left(\frac{D}{f_{\chi}}\right)^*} \not\in {\mathbb N},
 \]
and the corresponding residue vanishes in these cases as well.

\smallskip
%%%%%
Finally, for $\chi\in C_0$ and $m=(D/f_{\chi})^*$ we have
\[ 
n_{m^{\lambda_0}\alpha_{\nu}}(\chi^*) = \nu.
\]
Hence, recalling (\ref{eq:sumres=0}) and (\ref{eq:residue2}), we conclude that for all integers $(\nu,MD)=1$, we have
\begin{equation}\label{eq:sumchi}
a(\nu) \sum_{\chi\in C_0} \overline{\ell(\chi)} \chi^*(\nu)=0,
\end{equation}
where in view of \eqref{eq:Bm'}
\[
\begin{split}
 \ell(\chi) &= c(F^{\chi^*}) f(\chi, (D/f_{\chi})^*, 1/D) \frac{B_{(D/f_{\chi})^*}(F)}{((D/f_\chi)^*)^{s_0}} \\ &=
 c(F^{\chi^*}) f(\chi, (D/f_{\chi})^*, 1/D)  \prod_{p|(D/f_{\chi})} \frac{A_{\partial_p}(p)}{p^{\partial_p s_0}}.
 \end{split}
 \]
But (\ref{eq:Ap}) and Lemma \ref{lem:fchi} imply that $\ell(\chi)\neq 0$ for all $\chi\in C_0$, and since $F$ is saturated $(\bmod \, D)$, (\ref{eq:sumchi}) implies that the characters $\chi^*$ with $\chi\in C_0$ are linearly dependent, a contradiction. This shows that $d_{\chi^*}=2$ for all $\chi(\bmod \, D)$.

\smallskip
%%%%%
The proof that $\theta_{\chi^*}=\theta_F$ for all $\chi(\bmod \, D)$ can be carried out analogously. As before, we proceed by contradiction, keeping the notation introduced in the first part of the proof. Now $d_0=2$, $\lambda_0=1/2$, and $A_0$ contains all characters $\chi(\bmod \, D)$. We redefine $\theta_0$ to be any $\theta_{\chi^*}\neq\theta_F$, so that $F(s,1/D,\alpha,1/2)$ is holomorphic at $s=s_0$ by Lemma 4.
Hence we have that (\ref{eq:sumres=0}) still holds with $\lambda_0$ replaced by $1/2$. Then, repeating the subsequent calculations with $\lambda_0=1/2$, we arrive at (\ref{eq:sumchi}), which again gives a contradiction. The theorem now follows. \qed

\bigskip
\bigskip
\noindent
Jerzy Kaczorowski, Faculty of Mathematics and Computer Science, A.Mickiewicz University, 61-614 Pozna\'n, Poland. e-mail: \url{kjerzy@amu.edu.pl}

\medskip
\noindent
Alberto Perelli, Dipartimento di Matematica, Universit\`a di Genova, via Dodecaneso 35, 16146 Genova, Italy. e-mail: \url{perelli@dima.unige.it}

\end{document}